\documentclass[reqno,12pt,a4letter]{amsart}
\usepackage{amsmath, amsxtra, amssymb, latexsym, amscd, amsthm}
\usepackage{graphicx, color}
\usepackage[active]{srcltx}
\usepackage[utf8]{inputenc}
\usepackage[mathscr]{euscript}
\usepackage{mathrsfs,cite}
\usepackage[english]{babel}
\usepackage{color}
\usepackage[active]{srcltx}
\def\disp{\displaystyle}
\def\h{\hfill\triangle}

\newtheorem {theorem}{Theorem}[section]
\newtheorem {corollary}{Corollary}[section]

\newtheorem {lemma}{Lemma}[section]
\newtheorem {example}{Example}[section]
\newtheorem {definition}{Definition}[section]
\newtheorem {remark}{Remark}[section]
\def\ar{a\kern-.370em\raise.16ex\hbox{\char95\kern-0.53ex\char'47}\kern.05em}
\def\ees{{\accent"5E e}\kern-.385em\raise.2ex\hbox{\char'23}\kern-.08em}
\def\eex{{\accent"5E e}\kern-.470em\raise.3ex\hbox{\char'176}}
\def\AR{A\kern-.46em\raise.80ex\hbox{\char95\kern-0.53ex\char'47}\kern.13em}
\def\EES{{\accent"5E E}\kern-.5em\raise.8ex\hbox{\char'23 }}
\def\EEX{{\accent"5E E}\kern-.60em\raise.9ex\hbox{\char'176}\kern.1em}
\def\ow{o\kern-.42em\raise.82ex\hbox{
\vrule width .12em height .0ex depth .075ex \kern-0.16em \char'56}\kern-.07em}
\def\OW{O\kern-.460em\raise1.36ex\hbox{
\vrule width .13em height .0ex depth .075ex \kern-0.16em \char'56}\kern-.07em}
\def\UW{U\kern-.42em\raise1.36ex\hbox{
\vrule width .13em height .0ex depth .075ex \kern-0.16em \char'56}\kern-.07em}
\def\DD{D\kern-.7em\raise0.4ex\hbox{\char '55}\kern.33em}
\title{Existence of efficient and properly efficient solutions to problems of constrained vector optimization}

\author{DO SANG KIM$^1$}
\address{$^1$Department of Applied Mathematics, Pukyong National University, Busan 48513, Republic of Korea}
\email{dskim@pknu.ac.kr}

\author{BORIS S. MORDUKHOVICH$^{2,3}$}
\address{$^2$Department of Mathematics, Wayne State University, Detroit, MI 48202, USA}
\address{$^3$RUDN University, Moscow 117198, Russia}
\email{boris@math.wayne.edu}

\author{TI\EES N-S\OW N PH\d{A}M$^4$}
\address{$^4$Department of Mathematics, University of Dalat, 1 Phu Dong Thien Vuong, Dalat, Vietnam}
\email{sonpt@dlu.edu.vn}

\author{NGUYEN VAN TUYEN$^{5,6}$}
\address{$^5$School of Mathematical Sciences, University of Electronic Science and Technology of China, Chengdu, P.R. China}
\address{$^6$Department of Mathematics, Hanoi Pedagogical University 2, Xuan Hoa, Phuc Yen, Vinh Phuc, Vietnam}
\email{tuyensp2@yahoo.com; nguyenvantuyen83@hpu2.edu.vn}

\thanks{The first author was supported by the National Research Foundation of Korea Grant funded by the Korean Government (NRF-2016R1A2B4011589). Research of second author was partially supported by the USA National Science Foundation under grant DMS-1512846 and by the USA Air Force Office of Scientific Research under grant \#15RT0462. The third author was partially supported by Vietnam National Foundation for Science and Technology Development (NAFOSTED) under grant \#101.04-2016.05.}

\date{\today}



\begin{document}

\maketitle

\begin{abstract}
The paper is devoted to the existence of global optimal solutions for a general class of nonsmooth problems of constrained vector optimization without boundedness assumptions on constraint sets. The main attention is paid to the two major notions of optimality in vector problems: Pareto efficiency and proper efficiency in the sense of Geoffrion. Employing adequate tools of variational analysis and generalized differentiation, we first establish relationships between the notions of properness, $M$-tameness, and the Palais--Smale conditions formulated for the restriction of the vector cost mapping on the constraint set. These results are instrumental to derive verifiable necessary and sufficient conditions for the existence of Pareto efficient solutions in vector optimization. Furthermore, the developed approach allows us to obtain new sufficient conditions for the existence of Geoffrion-properly efficient solutions to such constrained vector problems.
\end{abstract}

\section{Introduction} \label{introduction}

This paper concerns some fundamental issues of global vector optimization that are revolved around the existence of efficient and properly efficient solutions under unbounded constraints. Such issues have been addressed in many publications; see, e.g., the books \cite{jahn04,Luc1989,Sawaragi1985} and the papers \cite{bao07,bao10,borwein83,Flores2002,Geoffrion68,Gutierrez2014,ha06,ha10} with the references therein. We offer here a new approach to these topics that allows us to derive significantly new existence theorems for a general class of problems in vector optimization. This approach is mainly based on advanced tools of variational analysis and generalized differentiation that provides essential improvements of known results even in the case of problems with smooth data.

The basic problem under consideration is formulated as follows:
\begin{equation}\label{problem}
\text{\rm Min}_{\,\mathbb{R}^m_+}\big\{f(x)\,\big|\,x\in\Omega\big\},\tag{VP}
\end{equation}
where $f\colon\mathbb{R}^n\to\mathbb{R}^m$ is a locally Lipschitz mapping, where $\Omega\subset\mathbb{R}^n$ is a nonempty and closed (not necessarily bounded) set, and where ``minimization" is understood in conventional terms of vector optimization that are specified below.

In the case of unconstrained problems \eqref{problem} with $\Omega=\mathbb{R}^n$, existence theorems for weak Pareto/weak efficient and the so-called relative Pareto (while not Pareto efficient) solutions were obtained in \cite{bao07,bao10,Mordukhovich2018} by using appropriate set-valued extensions of the Ekeland variational principle under the following major assumptions:
\begin{itemize}
\item[$\bullet$] $f$ is quasibounded from below, i.e., there exists a set $M\subset\mathbb{R}^m$ such that
\begin{eqnarray*}
f(\mathbb{R}^n)\subset M+\mathbb{R}^m_+.
\end{eqnarray*}
\item[$\bullet$] $f$ satisfies a certain Palais--Smale condition.
\end{itemize}
Somewhat related results for weak Pareto minimizers were obtained \cite{ha06} under more restrictive assumptions. As discussed in \cite{Kim2018}, such assumptions are rather limited. To improve them, powerful methods of semialgebraic geometry and polynomial optimization were invoked in \cite{Kim2018}. In this way the equivalence between the following conditions was proved therein when $\Omega=\mathbb{R}^n$ and $f$ is {\em polynomial} in \eqref{problem}; see below for the exact definitions:

\begin{itemize}
\item[$\bullet$] $f$ is proper at the sublevel $\bar y$.
\item[$\bullet$] $f$ satisfies the Palais--Smale condition at the sublevel $\bar{y}.$
\item[$\bullet$] $f$ satisfies the weak Palais--Smale condition at the sublevel $\bar{y}.$
\item[$\bullet$] $f$ is $M$-tame at the sublevel $\bar{y}.$
\end{itemize}
As consequences of these results, some sufficient conditions for the existence of Pareto efficient solutions of the unconstrained polynomial problem \eqref{problem} were given in \cite{Kim2018}.\vspace*{0.1in}

The main {\bf contributions} of this paper are significantly different from \cite{Kim2018}. First of all, we study the {\em constrained} problem \eqref{problem} with an arbitrary closed constraint set $\Omega$ and without any polynomial requirement on $f$, which is now replaced by {\em local Lipschitz continuity}. To proceed, we do not use methods of semialgebraic geometry but employ instead tools of {\em variational analysis and generalized differentiation}. Our {\em major results} are as follows:

{\bf(a)} Assuming that the image set $f(\Omega)$ has a {\em bounded section} at some $\bar{y}\in f(\Omega)$, which is indeed {\em necessary} for the existence of Pareto efficient solutions to \eqref{problem}, we show that the following statements are {\em equivalent}:
\begin{itemize}
\item the restriction $f|_\Omega$ of $f$ on $\Omega$ is proper at the sublevel $\bar y$.
\item the restriction $f|_\Omega$ satisfies the Palais--Smale condition at the sublevel $\bar{y}.$
\item the restriction $f|_\Omega$ satisfies the weak Palais--Smale condition at the sublevel $\bar{y}.$
\item the restriction $f|_\Omega$ is $M$-tame at the sublevel $\bar{y}.$
\end{itemize}

{\bf(b)} Based on these results, we obtain {\em necessary and sufficient} conditions for the existence of {\em Pareto efficient solutions} to problem~\eqref{problem}. As a byproduct of our approach, new sufficient conditions for the existence of {\em Geoffrion-properly efficient} solutions to \eqref{problem} are also derived.\vspace*{0.05in}

The rest of the paper is organized as follows. In Section~\ref{Section2} we recall some definitions and preliminary results from variational analysis and generalized differentiation. Section~\ref{Section3} is devoted to establishing relationships between properness, Palais--Smale conditions, and $M$-tameness. In Section~\ref{Section 4} we prove the existence of Pareto efficient and Geoffrion-properly efficient solutions to the vector optimization problem \eqref{problem}. The concluding Section~5 contains discussions of open problems to address in our future research.

\section{Preliminaries}\label{Section2}

Our notation is terminology are standard in variational analysis and vector optimization; see, e.g., the books \cite{jahn04,Mordukhovich2006,Rockafellar1998}. Recall that for any number $n\in\mathbb{N}:=\{1,2,\ldots\}$ we denote $x:=(x_1,\ldots,x_n)$ and equip the space $\mathbb{R}^n$ with the usual scalar product $\langle\cdot,\cdot\rangle$ and the Euclidean norm $\|\cdot\|$. The  closed unit ball in $\mathbb{R}^n$ is denoted by $\mathbb{B}^n$.

\subsection{Definitions of optimal solutions} Let
\begin{eqnarray*}
\mathbb{R}^m_+&:=&\big\{y:=(y_1,\ldots,y_m)\in\mathbb{R}^m\,\big|\,\,y_i\geq 0,\;i=1,\ldots,m\big\},\\
\mathrm{int}\,\mathbb{R}^m_+&:=&\big\{y:=(y_1,\ldots,y_m)\in\mathbb{R}^m\,\big|\,\,y_i>0,\;i=1,\ldots,m\big\}.
\end{eqnarray*}
Then $a\leqq b$ means $b-a\in\mathbb{R}^m_+$, $a\le b$ means $b-a\in\mathbb{R}^m_+\setminus\{0\}$, and $a<b$ means $b-a\in\text{int}\,\mathbb{R}^m_+$ for any vectors $a,b\in\mathbb{R}^m$.

\begin{definition}{\rm Given $\bar x\in\Omega$, we say that

\begin{enumerate}
\item [\bf(i)] $\bar{x}$ is a {\em Pareto efficient solution} to \eqref{problem} if there is no $x\in\Omega$ such that
\begin{eqnarray*}
f(x)&\leqq&f(\bar{x})\quad\textrm{ and }\quad f(x)\ne f(\bar{x}).
\end{eqnarray*}

\item [\bf(ii)] $\bar x$ is a {\em Geoffrion-properly efficient solution} to \eqref{problem} if it is a Pareto efficient solution and there
is a real number $M>0$ such that whenever $i\in\{1,\ldots,m\}$ and $x\in\Omega$ satisfying $f_i(x)<f_i(\bar x)$ there exists an index $j\in\{1,\ldots,m\}$ with $f_j(\bar x)<f_j(x)$ and
\begin{eqnarray*}
\frac{f_i(\bar x)-f_i(x)}{f_j(x)-f_j(\bar x)}&\le& M.
\end{eqnarray*}
\end{enumerate}
}\end{definition}

It follows from the definitions that every Geoffrion-properly efficient solution is a Pareto efficient solution to \eqref{problem} but not vice versa; see e.g., Example~\ref{Example41} below.

\subsection{Normals and subdifferentials}
Here we recall the notions of the normal cones to closed sets and the subdifferential of real-valued functions used in this paper. The reader is referred to \cite{Mordukhovich2006,Rockafellar1998} for more details.

\begin{definition}{\rm Consider a set $\Omega\subset\mathbb{R}^n$ and a point $\bar x\in\Omega.$
\begin{enumerate}
\item[\bf(i)] The {\em regular normal cone} (known also as the prenormal or Fr\'echet normal cone) $\widehat{N}(\bar x;\Omega)$ to
$\Omega$ at $\bar x$ consists of all vectors $v\in\mathbb{R}^n$ satisfying
\begin{eqnarray*}
\langle v,x-\bar{x}\rangle&\le& o(\|x-\bar{x}\|)\quad\textrm{ as }\quad x\to\bar{x}\quad\textrm{ with }\quad x\in\Omega.
\end{eqnarray*}

\item[\bf(ii)] The {\em limiting normal cone} (known also as the basic or Mordukhovich normal cone) $N(\bar x;\Omega)$ to $\Omega$ at
$\bar x$ consists of all vectors $v\in\mathbb{R}^n$ such that there are sequences $x^k\to\bar{x}$ with $x^k\in\Omega$ and $v^k\rightarrow v$ with $v^k\in\widehat N(x^k;\Omega)$ as $k\to\infty$.
\end{enumerate}}
\end{definition}

\begin{definition}{\rm Consider a function $\phi\colon\mathbb{R}^n\to\mathbb{R}$ and a point $\bar x\in\mathbb{R}^n$. The (limiting) {\em subdifferential} of $\phi$ at $\bar x$ is defined by
\begin{eqnarray*}
\partial\phi(\bar x)&:=&\big\{v\in\mathbb{R}^n\;\big|\;(v,-1)\in N\big((\bar x,\phi(\bar x));\mathrm{epi}\,\phi\big)\big\}
\end{eqnarray*}
via the limiting normal cone to the epigraph $\mathrm{epi}\,\phi$ of $\phi$ given by
\begin{eqnarray*}
\mathrm{epi}\,\phi&:=&\big\{(x,y)\in\mathbb{R}^n\times\mathbb{R}\,\big|\,\phi(x)\le y\big\}.
\end{eqnarray*}}
\end{definition}

In \cite{Mordukhovich2006,Mordukhovich2018,Rockafellar1998} the reader can find equivalent analytic descriptions of the subdifferential $\partial\phi(\bar x)$ and comprehensive studies of it and related constructions. In the case of convex sets and functions the above normal cone and subdifferential notions reduce to the corresponding concepts of convex analysis. Furthermore, we have $\partial\phi(\bar x)=\{\nabla\phi(\bar x)\}$ if $\phi$ is strictly differentiable at $\bar x$; in particular, when it is smooth around this point.\vspace*{0.05in}

Next we present several known statements, which play significant roles in the proofs of the main results. The first lemma is the classical subdifferential formula of convex analysis.

\begin{lemma}\label{Lemma21}
For each $\bar{x}\in\mathbb{R}^n$ we have
\begin{eqnarray*}
\partial\left(\|\,\cdot\,-\bar{x}\|\right)(x)&=&
\begin{cases}
\disp\frac{x-\bar{x}}{\|x-\bar{x}\|}&\textrm{ if }\;x\ne\bar{x},\\
\mathbb{B}^n&\text{ otherwise.}
\end{cases}
\end{eqnarray*}
\end{lemma}

The following major results of subdifferential calculus and necessary optimality conditions for scalar nonsmooth optimization are used below in our derivation of the existence theorems for \eqref{problem} even in the case of problems with smooth initial data.

\begin{lemma}[{see \cite[Theorem~3.36]{Mordukhovich2006}}]\label{Lemma22}
Let the functions $\phi_i\colon\mathbb{R}^n \rightarrow{\mathbb{R}}$, $i=1,\ldots,m$, be locally Lipschitzian around $\bar{x}\in\mathbb{R}^n$. Then we have the subdifferential sum rule
\begin{eqnarray*}
\partial(\phi_1+\cdots+\phi_m)(\bar x)&\subset&\partial\phi_1(\bar{x})+\cdots+\partial\phi_m(\bar{x}).
\end{eqnarray*}
\end{lemma}

\begin{lemma}[{see \cite[Theorem~3.46]{Mordukhovich2006}}]\label{Lemma23}
Let $\phi_1,\ldots,\phi_m\colon{\Bbb R}^n\rightarrow{\Bbb R}$ be locally Lipschitzian around $\bar{x}\in\mathbb{R}^n$. Then the maximum function
\begin{eqnarray*}
\phi(x):=\max_{1\le i\le m}\phi_i(x),\quad x\in\mathbb{R}^n,
\end{eqnarray*}
is locally Lipschitzian around $\bar{x}$, and we have the inclusion
\begin{eqnarray*}
\partial\phi(\bar{x})&\subset& \left\{\disp\sum_{i\in I(\bar{x})}\alpha_i\partial\phi_i(\bar{x})\;\Big|\;\alpha_i\ge 0\quad\textrm{ and }\quad\disp\sum_{i\in I(\bar{x})}\alpha_i=1\right\},
\end{eqnarray*}
where the active index set is defined by $I(\bar{x}):=\{i\;|\;\phi_i(\bar{x})=\phi(\bar{x})\}$.
\end{lemma}

\begin{lemma}[{see\cite[Theorem~5.21(iii)]{Mordukhovich2006}}]\label{Lemma24}
Let the function $\phi_i\colon\mathbb{R}^n\rightarrow{\mathbb{R}}$, $i=0,\ldots,m$, be locally Lipschitzian $\bar x\in\Omega$, and let the set $\Omega\subset\mathbb{R}^n$ be locally closed around this point. If $\bar x$ is a local minimizer of the function $\phi_0$ on the set
\begin{eqnarray*}
\big\{x\in\Omega\,\big|\,\phi_i(x)\le 0,\;i=1,\ldots,m\big\},
\end{eqnarray*}
then there exist nonnegative numbers $\lambda_i,\;i=0,\ldots,m$, with $\sum_{i=1}^m\lambda_i=1$ such that
\begin{eqnarray*}
&&0\in\disp\sum_{i=0}^{m}\lambda_i\partial\phi_i(\bar x) + N(\bar x;\Omega),\\
&&\lambda_i\phi_i(\bar x)=0,\;i=1,\ldots,m.
\end{eqnarray*}
\end{lemma}

\section{Properness, Palais--Smale conditions, and M-Tameness}\label{Section3}

Let $f:=(f_1,\ldots, f_m)\colon\mathbb{R}^n\to\mathbb{R}^m$ be a locally Lipschitz mapping, and let $\Omega\subset\mathbb{R}^n$ be a nonempty and  closed set. In this section we establish close relationships between the properness, Palais--Smale conditions, and $M$-tameness for the restriction $f|_{\Omega}$ of $f$ on $\Omega$. To proceed, the following notion is needed.

\begin{definition}{\rm
Let $A$ be a subset in $\mathbb{R}^m$, and let $y\in\mathbb{R}^m$. The set $A\cap (y-\mathbb{R}^m_+)$ is called a {\em section} of $A$ at $y$ and is denoted by $[A]_y.$ The section $[A]_y$ is said to be {\em bounded} if there exists a vector $a\in\mathbb{R}^m$ such that
\begin{eqnarray*}
[A]_y&\subset&a+\mathbb{R}^m_+.
\end{eqnarray*}
}\end{definition}

It is easy to observe that for a Pareto efficient solution $\bar x$ to \eqref{problem} we have
\begin{eqnarray*}
[f(\Omega)]_{f(\bar x)}&=&\big\{f(\bar x)\big\}.
\end{eqnarray*}
Thus the condition that $f(\Omega)$ admits at least one bounded section is {\em necessary for the existence of Pareto efficient solutions} to \eqref{problem}.\vspace*{0.05in}

Next we introduce the notions of {\em properness} for the restricted cost mapping, which are instrumental to prove the existence of optimal solutions to \eqref{problem}.

\begin{definition}{\rm We say that:
\begin{enumerate}
\item[\bf(i)] The restriction $f|_\Omega$ of $f$ on $\Omega$ is {\em proper at sublevel} $y\in\mathbb{R}^m$ if
\begin{equation*}
\forall\,\{x^k\}\;\subset\;\Omega,\;\|x^k\|\to\infty,\;f(x^k)\;\leqq\;y\Longrightarrow\|f(x^k)\|\to\infty\;\mbox{ as }\;k\to\infty.
\end{equation*}

\item[\bf(ii)] The restriction $f|_{\Omega}$  is {\em proper} if it is proper at every sublevel $y\in\mathbb{R}^m$.
\end{enumerate}
}\end{definition}

For each $\bar y\in\big(\mathbb{R}\cup\{\infty\}\big)^m$, consider the sets
\begin{eqnarray*}
\widetilde{K}_{\infty,\leqq \bar{y}}\,(f,\Omega)&:=&\big\{y\in\mathbb{R}^m \big|\,\exists\,\{x^k\}\subset\Omega,\;f(x^k)\leqq\bar y,\;\|x^k\|\to\infty, f(x^k)\to y,\;\mbox{ and}\\
&&\qquad\qquad\qquad\nu(x^k)\to 0\;\mbox{ as }\,k\to\infty\},\\
{K}_{\infty,\leqq\bar{y}}\,(f,\Omega)&:=&\big\{y\in\mathbb{R}^m \,\big|\,\exists\,\{x^k\}\subset\Omega,\;f(x^k)\leqq\bar y,\;\|x^k\|\to\infty,\;f(x^k)\to y,\;\mbox{ and}\\
&&\qquad\qquad\qquad\|x^k\|\nu(x^k)\to 0\;\mbox{ as }\;k\to\infty\},
\end{eqnarray*}
where $\nu\colon\mathbb{R}^n\to\mathbb{R}$ is the (extended) {\em Rabier function} defined by
\begin{eqnarray}\label{Eqn1}
\nu(x):=\inf\Big\{\Big\|\sum_{i=1}^m\lambda_i v^i+w\Big\|\;\Big|\;v^i\in\partial f_i(x)\cup\partial(-f_i)(x),\,w\in N(x;\Omega),\,\lambda\in\mathbb{R}_+^m,\,\disp\sum_{i=1}^{m}\lambda_i=1\Big\}.
\end{eqnarray}

Following \cite[Chapter~2]{HHVui2017}, we also consider the set
\begin{eqnarray*}
T_{\infty,\leqq\bar y}\,(f,\Omega)&:=&\big\{y\in\mathbb{R}^m\;\big|\;\exists\,\{x^k\}\subset\Gamma(f,\Omega),\;f(x^k)\leqq\bar y,\;\|x^k\|\to\infty,\; \mbox{ and }\\
&&\qquad\qquad\qquad f(x^k)\to y\;\mbox{ as }\;k\to\infty\},
\end{eqnarray*}
defined via the construction
\begin{eqnarray*}
\Gamma(f,\Omega)\ :=\ \Big\{x\in\Omega&\Big|&\exists v^i\in\partial f_i(x)\cup\partial(-f_i)(x) \ \textrm{ for } \ i=1,\ldots,m,\\
&&\exists(\lambda,\mu)\in\mathbb{R}_+^m\times\mathbb{R}\ \textrm{ with } \ \disp\sum_{i=1}^m\lambda_i+|\mu|=1 \textrm{ such that }\\
&&0\ \in\disp\sum_{i=1}^m\lambda_i v^i+\mu x+N(x;\Omega)\Big\}.
\end{eqnarray*}
\begin{remark}\label{rem}{\rm 
It follows from the proofs given below that the results of Theorems~\ref{Theorem31}, \ref{Theorem42} and Corollaries~\ref{Corollary31}, \ref{Corollary41} hold true if in the definitions of the function $\nu$ and the set $\Gamma(f,\Omega)$ we delete the terms $\partial(-f_i)(x)$. We include this in the original definitions for simplicity and for the unification with other results of the paper.}
\end{remark}

When $\bar{y}=(\infty,\ldots,\infty)$, we simplify the notation by writing $\widetilde{K}_{\infty}(f,\Omega),$ ${K}_{\infty}(f,\Omega),$ and $T_\infty(f,\Omega)$ instead of $\widetilde{K}_{\infty,\leqq \bar{y}}(f,\Omega),$ $K_{\infty,\leqq\bar{y}}(f,\Omega),$ and $T_{\infty,\leqq\bar{y}}(f,\Omega)$, respectively.

It follows from the definitions that the properness of $f|_{\Omega}$ at sublevel $\bar y\in\mathbb{R}^m$ yields
\begin{eqnarray*}
T_{\infty,\leqq\bar{y}}(f,\Omega)&=&\widetilde{K}_{\infty,\leqq\bar{y}}\,(f,\Omega)\ = \ {K}_{\infty,\leqq\bar{y}}\,(f,\Omega)\ =\ \emptyset.
\end{eqnarray*}
The {\em converse does not hold} in general. Indeed, let $\Omega:=\mathbb{R}^2$, and let $f\colon\mathbb{R}^2\to\mathbb{R}$ be a real-valued function defined by $f(x_1,x_2):=x_1+x_2$. It is easy to check that
\begin{eqnarray*}
T_{\infty,\leqq \bar y}\,(f,\Omega)=\widetilde{K}_{\infty,\leqq \bar{y}}\,(f,\Omega)={K}_{\infty,\leqq\bar{y}}\,(f,\Omega)=\emptyset\quad\textrm{ for all }\quad \bar y\in\mathbb{R},
\end{eqnarray*}
while $f$ is not proper at every sublevel. Nevertheless, we have the following rather surprising result the proof of which is based on {\em variational arguments} and {\em subdifferential calculus}.

\begin{theorem}\label{Theorem31}
Assume that there exists $\bar y\in f(\Omega)$ such that the section $[f(\Omega)]_{\bar y}$ is bounded.
Then  the following statements are equivalent:
\begin{enumerate}
\item[\bf(i)] $f|_{\Omega}$ is proper at the sublevel $\bar y$.
\item[\bf(ii)] $f|_{\Omega}$ satisfies the {\sc Palais--Smale condition} at $\bar y$,
i.e., $\widetilde{K}_{\infty,\leqq\bar{y}}\,(f,\Omega)=\emptyset$.
\item[\bf(iii)] $f|_{\Omega}$ satisfies the {\sc weak Palais--Smale condition} at $\bar y$, i.e.,
${K}_{\infty,\leqq\bar{y}}\,(f,\Omega)=\emptyset$.
\item [\bf(iv)] $f|_{\Omega}$ is {\sc M-tame} at the sublevel $\bar{y}$, i.e., $T_{\infty,\leqq\bar{y}}\,(f,\Omega)=\emptyset$.
\end{enumerate}
Furthermore, the set $[f(\Omega)]_{\bar y}$ is nonempty and compact provided that one of the above equivalent conditions is satisfied.
\end{theorem}
{\bf Proof.} Note first that implications (i)$\Rightarrow$(ii), (ii)$\Rightarrow$(iii), and (i)$\Rightarrow$(iv) are obvious.

To prove (iii)$\Rightarrow$(i), we argue by contradiction and assume that $f$ is not proper at the sublevel $\bar y$. Since the section
$Y:=[f(\Omega)]_{\bar y}$ is bounded, there exists $\widehat y\in\mathbb{R}^m$ with $\widehat y\not\in Y$ and $\widehat y\leqq\bar{y}$. Consider the nonempty set
\begin{eqnarray}\label{Eqn2}
X&:=&f^{-1}(Y)\cap\Omega=\big\{x\in\Omega\,\big|\,f(x)\leqq\bar{y}\big\}\ne\emptyset
\end{eqnarray}
and define the function $\phi\colon\mathbb{R}^n\rightarrow\mathbb{R}$ by
\begin{eqnarray}\label{Eqn3}
\phi(x)&:=&\max_{i=1,\ldots,m}|f_i(x)-\widehat y_i|.
\end{eqnarray}
It is clear that the function $\phi$ is nonnegative, locally Lipschitzian and satisfies the condition $\phi(x)>0$ for all $x\in X$. Since $f$ is not proper at the sublevel $\bar y$, we see that the number
\begin{eqnarray}\label{Eqn4}
c&:=&\liminf\limits_{x\in X,\,\|x\|\to\infty}\phi(x)
\end{eqnarray}
is finite. For each $R>0$ consider the quantity
\begin{eqnarray*}
\frak{m}(R)&:=&\inf\limits_{x\in X,\,\|x\|\ge R}\phi(x)
\end{eqnarray*}
and observe that $\frak{m}$ is a nondecreasing and nonnegative function with $\lim\limits_{R\to\infty}\frak{m}(R)=c.$ Thus for each $k\in\mathbb{N}$ there exists $R_k>k$ satisfying
\begin{eqnarray*}
\frak{m}(R)&\ge& c-\frac{1}{2k}\quad\textrm{whenever}\quad R\ge R_k.
\end{eqnarray*}
Choose now $x^k\in X$ with $\|x^k\|>2R_k$ and such that
\begin{eqnarray*}
\phi (x^k)&<&\frak{m}(2 R_k)+\frac{1}{2k}.
\end{eqnarray*}
We clearly have the chain of inequalities
\begin{eqnarray*}
\frak{m}(2R_k) \ \le \ \phi(x^k)&<&\frak{m}(2 R_k)+\frac{1}{2k}\\
&\le&c+\frac{1}{2k} \\
&\le&\frak{m}(R_k)+\frac{1}{k}\\
&=&\inf\limits_{x\in X,\,\|x\|\ge R_k}\phi(x)+\frac{1}{k}.
\end{eqnarray*}
We are now in a position to apply the Ekeland variational principle \cite{Ekeland1979} (see, e.g.,  \cite[Theorem~2.26]{Mordukhovich2006}) to the function $\phi$ on the set $\{x\in X\,|\,\|x\|\ge R_k\}$ with the parameters $\varepsilon:=\frac{1}{3k}$ and $\lambda:=\dfrac{\|x^k\|}{2}$ therein. Note that in the finite-dimensional setting under consideration this result and other variational principles can be proved easily; see \cite[Theorem~2.12]{Mordukhovich2018}. In this way we find $u^k\in X$ with $\|u^k\|\ge R_k$ satisfying the following conditions:
\begin{itemize}
\item [(a)] $\frak{m}(R_k)\le\phi(u^k)\leq\phi(x^k),$
\item [(b)] $\|u^k - x^k\|\le\lambda$, and
\item [(c)] $\phi (x)+\dfrac{\varepsilon}{\lambda}\|x-u^k\|\ge\phi(u^k)$ for all $x\in X$ with $\|x\|\ge R_k.$
\end{itemize}
It follows from (a) that $\phi(u^k)\to c$ as $k\to\infty$, while (b) yields
\begin{eqnarray*}
R_k &<& \disp \dfrac{\|x^k\|}{2}\ \le\ \|u^k\| \ \le \ \dfrac{3}{2}\|x^k\|
\end{eqnarray*}
and implies, in particular, that $\|u^k\|\to\infty$ as $k\to\infty$. Applying the necessary optimality conditions from  Lemma~\ref{Lemma24} to the nonsmooth scalar optimization problem in (c) allows us to find $(\kappa,\beta)\in\mathbb{R}_+\times\mathbb{R}^m_+$ for which
\begin{eqnarray} \label{Eqn5}
&& 0 \in \kappa\,\partial\left[\phi(\cdot)+\dfrac{\varepsilon}{\lambda}\|\cdot-u^k\|\right](u^k)+\disp\sum_{i=1}^{m}\beta_i\partial f_i(u^k) + N(u^k;\Omega),\\
&& \beta_i \big (f_i(u^k) - \bar{y}_i \big) = 0 \ \textrm{ for } \ i = 1, \ldots, m, \quad  \textrm{ and } \quad 1 = \kappa+\disp\sum_{i=1}^{m}\beta_i. \nonumber
\end{eqnarray}
The subdifferential sum rule from Lemma~\ref{Lemma22} and the calculation of Lemma~\ref{Lemma21} give us
\begin{eqnarray} \label{Eqn6}
\partial\left[\phi(\cdot)+\disp\dfrac{\varepsilon}{\lambda}\|\cdot-u^k\|\right](u^k)&\subset&\partial\phi(u^k)+\dfrac{\varepsilon}{\lambda}\mathbb{B}^n.
\end{eqnarray}
In order to evaluate the subdifferential $\partial\phi(u^k)$ of the maximum function \eqref{Eqn3} in our setting, define the index sets
\begin{eqnarray*}
I_1^{\pm}&:=&\big\{i\,\big|\,f_i(u^k) - \widehat{y}_i = \pm\phi(u^k)\big\},\\
I_2&:=&\big\{i\,\big|\,f_i(u^k)-\bar{y}_i=0\big\}.
\end{eqnarray*}
Since $\phi(u^k)>0$, we get $I_1^-\cap I_1^+=\emptyset$. Furthermore, the condition $\widehat y\le\bar{y}$ ensures that $I_1^-\cap I_2=\emptyset$.
Then Lemma~\ref{Lemma23} tells us that
\begin{eqnarray*}
\partial\phi(u^k)&\subset&\left\{\disp\sum_{i\in I_1^+}\alpha_i\partial f_i(u^k) + \sum_{i\in I_1^-}\alpha_i \partial (-f_i)(u^k)\,\Big|\,\alpha_i\ge 0,\; \disp\sum_{i\in I_1^+\cup I_1^-}\alpha_i=1\right\}.
\end{eqnarray*}
This, together with \eqref{Eqn5} and \eqref{Eqn6}, implies that there exist numbers $\alpha_i\ge 0$ for $i\in I_1^+\cup I_1^-$ with $\sum_{i\in I_1^+\cup I_1^-}\alpha_i=1$ such that
\begin{eqnarray*}
0&\in&\disp\sum_{i\in I_1^+}\kappa\alpha_i\partial f_i(u^k) + \sum_{i\in I_1^-}\kappa\alpha_i\partial (-f_i)(u^k)+\sum_{i=1}^{m}\beta_i\partial f_i(u^k) + N(u^k;\Omega)+\kappa\dfrac{\varepsilon}{\lambda}\mathbb{B}^n.
\end{eqnarray*}
For each $i\in I^+_1\cup I^-_1\cup I_2$, we denote
\begin{eqnarray*}
\lambda_i&:=&
\begin{cases}
\kappa\alpha_i&\textrm{ if }\;i\in I^+_1\setminus I_2,\\
\kappa\alpha_i+\beta_i&\textrm{ if }\;i\in I^+_1\cap I_2,\\
\beta_i&\textrm{ if }\;i\in I_2\setminus I^+_1,\\
\kappa\alpha_i &\textrm{ if }\;i\in I^-_1,\\
0&\textrm{ otherwise,}
\end{cases}
\end{eqnarray*}
and arrive at the following relationships
\begin{eqnarray*}
&&0 \ \in \ \disp \sum_{i \not\in I_1^{-}}\lambda_i\partial f_i(u^k)+\sum_{i \in I_1^{-}}\lambda_i\partial(-f_i)(u^k)+ N(u^k;\Omega)+\kappa\dfrac{\varepsilon}{\lambda}\mathbb{B}^n,\\
&&\lambda_i\ge 0 \textrm{ for } i=1,\ldots,m,\quad \textrm {and }\quad\disp\sum_{i=1}^m\lambda_i  \ = \ \sum_{i\in I^+_1\cup I^-_1}\kappa\alpha_i+\sum_{i\in I_2} \beta_i \ = \ \kappa+\disp\sum_{i\in I_2}\beta_i  \ = \ 1.
\end{eqnarray*}

Observe that if we choose $\widehat y$ so that $\widehat y < y$ for all $y\in Y,$ then $I_1^{-}=\emptyset;$ cf.\ Remark~\ref{rem}. It follows from the definition of the Rabier function $\nu$ that
\begin{eqnarray*}
\nu(u^k) &\le& \kappa\frac{\varepsilon}{\lambda} \ \le \ \frac{\varepsilon}{\lambda}  \ = \ \frac{2}{3 k \|x^k\|} \ \le \ \frac{1}{k \|u^k\|}.
\end{eqnarray*}
Consequently, we get the estimate
\begin{eqnarray*}
\|u^k\|\nu(u^k)&\le&\disp\frac{1}{k}\;\mbox{ for each }\;k\in\mathbb{N},
\end{eqnarray*}
and therefore $\|u^k\|\nu(u^k)\to 0$ as $k\to\infty.$

On the other hand, it follows from the boundedness of the section $[f(\Omega)]_{\bar y}$ and the inclusion $\{f(u^k)\}\subset[f(\Omega)]_{\bar y}$ that the sequence $\{f(u^k)\}$ has an accumulation point, say $y\in\mathbb{R}^m$. Thus $y\in K_{\infty,\leqq\bar y}\,(f,\Omega),$ a contradiction that verifies implication (iii)$\Rightarrow$(i).

Next we prove (iv)$\Rightarrow$(i). Assume on the contrary that $f$ is not proper at the sublevel $\bar y$. Then there exists a sequence $\{x^k\}\subset\Omega$ such that $\|x^k\|\to\infty$ as $k\to\infty,$ $f(x^k)\leqq\bar y,$ and the sequence of images $\{f(x^k)\}$ is bounded.

Since the section $Y:=[f(\Omega)]_{\bar y}$ is bounded, there exists $\widehat y\in\mathbb{R}^m$ with $\widehat y\not\in Y$ and $\widehat y\le\bar{y}.$ As above, consider the set $X$ from \eqref{Eqn2}, the maximum function $\phi(x)$ from \eqref{Eqn3}, and then conclude that the number $c$ defined in \eqref{Eqn4} is finite.
For each $k\in\mathbb{N}$ we form the following scalar nonsmooth optimization problem:
\begin{eqnarray*}
&\textrm{minimize  }&\phi(x)\\
&\textrm{subject to }&x\in X\quad{\rm and }\quad\|x\|^2-\|x^k\|^2=0.
\end{eqnarray*}
Since the constraint set here is nonempty and compact, this problem admits an optimal solution denoted by $v^k$. The usage of necessary optimality conditions from Lemma~\ref{Lemma24} give us a triple $(\kappa,\beta,\mu)\in\mathbb{R}_+\times \mathbb{R}^m_+\times\mathbb{R}$ satisfying the relationships
\begin{eqnarray*}
&&0\in\kappa\,\partial\phi(v^k)+\disp\sum_{i=1}^m\beta_i\partial\big(f_i(\cdot)-\bar y_i \big)(v^k)+2\mu v^k+N(v^k;\Omega),\\
&&\beta_i\big(f_i(v^k)-\bar y_i\big)=0 \ \textrm{ for }\;i=1,\ldots,m\quad \textrm{ and }\quad 1=\kappa+\disp\sum_{i=1}^m\beta_i+|\mu|.
\end{eqnarray*}
Proceeding now as in the proof of the previous implication which taking into account the subdifferential sum rule and the subdifferential calculation for the maximum function together with the modified form of the cost function, we arrive at the conditions
\begin{eqnarray*}
&&0\ \in \ \disp\sum_{i\not\in I_1^{-}}\lambda_i\partial f_i(v^k)+\sum_{i\in I_1^{-}}\lambda_i\partial(-f_i)(v^k)+2\mu v^k+N(v^k;\Omega),\\
&&\disp\sum_{i=1}^m\lambda_i+2|\mu| \ = \ \sum_{i\in I^+_1\cup I^-_1}\kappa\alpha_i+\sum_{i\in I_2}\beta_i+2|\mu| \ = \ \kappa+\sum_{i\in I_2}\beta_i+2|\mu| \ = \ 1+|\mu| \ \ge \ 1
\end{eqnarray*}
with the same index sets and the expressions for $\lambda_i$ as above. We clearly get $v^k\in\Gamma(f,\Omega)$.

Thus we constructed the sequence $\{v^k\}$ with the following properties:
\begin{enumerate}
\item [(a)] $\{v^k\}\ \subset\ \Gamma(f,\Omega)$;
\item[(b)] $\|v^k\|\ =\ \|x^k\|\to\infty$ as $k\to\infty$;
\item[(c)] $\phi(v^k)\ \le \ \phi(x^k)$ for all $k\in\mathbb{N}$;
\item[(d)] $f(v^k)\ \leqq \ \bar y$ for all $k\in\mathbb{N}$.
\end{enumerate}
It follows from the boundedness of the section $[f(\Omega)]_{\bar y}$ and the inclusion $\{f(v^k)\}\subset[f(\Omega)]_{\bar y}$ that the sequence $\{f(v^k)\}$ has an accumulation point $y\in\mathbb{R}^m$. Therefore $y\in T_{\infty,\leqq\bar y}\,(f,\Omega)$, a contradiction. This completes the proof of the equivalence between all the properties (i)--(iv).

Let us finally verify the last statement of the theorem. Suppose that (i) holds and then show that the set $[f(\Omega)]_{\bar y}$ is closed and hence it is compact. To proceed, take an arbitrary sequence $\{y^k\}\subset[f(\Omega)]_{\bar y}$ converging to $y\in\mathbb{R}^m$ and find a sequence $\{x^k\}\subset\Omega$ such that $f(x^k)=y^k\leqq\bar{y}$ for all $k\in\mathbb{N}$. Since $\displaystyle\lim_{k\to\infty}f(x^k)=\displaystyle\lim_{k\to\infty}y^k=y$, it follows from (i) that the sequence $\{x^k\}$ is  bounded. Thus $\{x^k\}$ has an accumulation point $x$, which belongs to $\Omega$ due to the closedness of this set. The continuity of $f$ implies that $y=f(x)$, and  consequently we have that $y\in f(\Omega)$. Noting that $y\leqq{\bar y}$ gives us the inclusion $y\in[f(\Omega)]_{\bar y}$, which therefore completes the proof of the theorem. $\h$\vspace*{0.05in}

The results of Theorem~\ref{Theorem31} significantly extend the recent ones from \cite{Kim2018}, where such an equivalence is established in the case of $\Omega=\mathbb{R}^n$ and polynomial mappings $f$ by using methods of semialgebraic geometry. The proof of \cite{Kim2018} is based on the inclusion
\begin{eqnarray}\label{Eqn7}
T_{\infty, \bar y}(f,\mathbb{R}^n)&\subset&K_{\infty,\bar y}(f,\mathbb{R}^n)
\end{eqnarray}
valid when $f$ is polynomial. The following example shows that if $f$ is not polynomial, then \eqref{Eqn7} fails. Thus the approach of \cite{Kim2018} cannot be applied to our general setting, while the new approach of variational analysis allows to treat \eqref{problem} in full generality.

\begin{example}{\rm
Let $\bar y=0$, and let $f\colon\mathbb{R}\to\mathbb{R}$ be defined by $f(x):=\sin x$. We claim that $0\in T_{\infty,\bar y}(f,\mathbb{R}^n)\setminus K_{\infty,\bar y}(f,\mathbb{R}^n)$ and so $T_{\infty,\bar y}(f,\mathbb{R}^n)\nsubseteq K_{\infty,\bar y}(f,\mathbb{R}^n).$ Indeed, let $x^k=2k\pi$ for all $k\in\mathbb{N}$. It is easily seen that $\Gamma (f,\mathbb{R})=\mathbb{R}$ and hence $\{x^k\}\subset\Gamma (f,\mathbb{R})$. Since $x^k\to\infty$ and $f(x^k)\to 0$ as $k\to\infty$, we have that $0\in T_{\infty,\bar y}(f,\mathbb{R}^n)$. To show now that $0\notin K_{\infty,\bar y}(f,\mathbb{R}^n)$, assume the contrary and then find a sequence $\{u^k\}\subset\mathbb{R}$ such that $f(u^k)\leqq 0$, $u^k\to\infty$,  $f(u^k)\to 0$, and $u^k \nabla f(u^k)\to 0$ as $k\to\infty$. This implies that $\sin u^k\to 0$ and $u^k\cos u^k\to 0$ as $k\to\infty$. Using $u^k\to\infty$ and $u^k\cos u^k\to 0$, we get that $\cos u^k\to 0$ as $k\to\infty$. Consequently, $\sin^2 u^k+\cos^2 u^k\to 0$ as $k\to\infty$, which contradicts the fact that $\sin^2 u^k + \cos^2 u^k=1$ for all $k\in\mathbb{N}$ and thus verifies the failure of inclusion \eqref{Eqn7}.}
\end{example}

We conclude this section with an immediate consequence of Theorem~\ref{Theorem31}; cf.\ \cite[Theorem~2.5]{HHVui2017} for the case where $f$ is polynomial and $\Omega=\mathbb{R}^m$.

\begin{corollary} \label{Corollary31}
Assume that every section of the set $f(\Omega)$ is bounded. Then the following assertions are equivalent:
\begin{enumerate}
\item [\bf(i)] $f|_{\Omega}$ is proper.
\item[\bf(ii)] $f|_{\Omega}$ satisfies the Palais--Smale condition: $\widetilde{K}_{\infty}\,(f,\Omega)=\emptyset$.
\item [\bf (iii)] $f|_{\Omega}$ satisfies the weak Palais--Smale condition: ${K}_{\infty}\,(f,\Omega)=\emptyset$.
\item [\bf(iv)] $f|_{\Omega}$ is M-tame: $T_{\infty}\,(f,\Omega)=\emptyset$.
\end{enumerate}
Furthermore, every section of the set $f(\Omega)$ is compact provided that one of the above equivalent conditions is satisfied.
\end{corollary}

\section{Existence of optimal solutions}\label{Section 4}

This section contains our main results on the existence of optimal solutions to constrained vector optimization problem \eqref{problem} in the general nonsmooth setting. We start with deriving verifiable {\em necessary and sufficient} conditions for the existence of Pareto efficient solutions.

\subsection{Existence of Pareto efficient solutions}

Given $\bar y \in (\mathbb{R}\cup\{\infty\})^m$, denote
\begin{eqnarray*}
{K}_{0,\leqq\bar{y}}\,(f,\Omega)&:=&\big\{f(x)\in\mathbb{R}^m \,\big|\,x\in\Omega,\;f(x)\leqq\bar y,\;\nu(x)=0\big\},
\end{eqnarray*}
where $\nu(x)$ is the Rabier function defined in \eqref{Eqn1}. The motivation behind this definition comes from the observation that if $\bar x$ is a Pareto efficient solution to problem \eqref{problem}, then $\nu(\bar x)=0$ and so $f(\bar x)\in{K}_{0,\leqq\bar{y}}\,(f,\Omega)$ with
$\bar{y}=f(\bar x)$.

\begin{theorem}\label{Theorem41}
The following assertions are equivalent:
\begin{enumerate}
\item [\bf(i)] Problem \eqref{problem} admits a Pareto efficient solution.
\item [\bf(ii)] There exists a vector $\bar y\in f(\Omega)$ such that the section $[f(\Omega)]_{\bar y}$ is bounded
and the inclusion $\widetilde{K}_{\infty,\leqq\bar{y}}\,(f,\Omega)\subset {K}_{0,\leqq \bar{y}}\,(f,\Omega)$ holds.
\item [\bf(iii)] There exists a vector $\bar y\in f(\Omega)$ such that the section $[f(\Omega)]_{\bar y}$ is bounded and the
inclusion ${K}_{\infty,\leqq\bar{y}}\,(f,\Omega)\subset{K}_{0,\leqq  \bar{y}}\,(f,\Omega)$ holds.
\item [\bf(iv)] There exists a vector $\bar y\in f(\Omega)$ such that the section $[f(\Omega)]_{\bar y}$ is bounded and the
inclusion $T_{\infty,\leqq\bar{y}}\,(f,\Omega)\subset{K}_{0,\leqq \bar{y}}\,(f,\Omega)$ holds.
\end{enumerate}
\end{theorem}
{\bf Proof}. First we justify in parallel implications (i)$\Rightarrow$(ii), (i)$\Rightarrow$(iii), and (i)$\Rightarrow$ iv).  To this end, let $\bar{x}\in \Omega$ be a Pareto efficient solution to \eqref{problem}, and let $\bar{y}:= f(\bar{x})$.  As mentioned above, the section $[f(\Omega)]_{\bar y}$ is just $\{\bar{y}\}$ while containing in this case the sets ${K}_{0,\leqq\bar{y}}\,(f,\Omega),$  $\widetilde{K}_{\infty,\leqq\bar{y}}\,(f,\Omega),$ ${K}_{\infty,\leqq \bar{y}}\,(f,\Omega),$ and ${T}_{\infty,\leqq\bar{y}}\,(f,\Omega).$ Furthermore, the necessary optimality conditions from Lemma~\ref{Lemma24} ensure the existence of $\lambda\in\mathbb{R}^m_+$ such that $\sum_{i=1}^{m}\lambda_i=1$ and
\begin{eqnarray*}
0&\in&\disp\sum_{i=1}^{m}\lambda_i\partial f_i(\bar{x}) + N(\bar{x};\Omega).
\end{eqnarray*}

By definition \eqref{Eqn1} we have $\nu(\bar{x})=0$, and therefore $\bar{y}\in {K}_{0,\leqq\bar{y}}\,(f,\Omega)$. Thus the conditions in (ii), (iii), and (iv) follow immediately from these facts giving us necessary conditions for the existence of Pareto efficient solutions to \eqref{problem}.

Next we verify implications (ii)$\Rightarrow$(i), (iii)$\Rightarrow$(i), and (iv)$\Rightarrow$(i), which justify the sufficiency of conditions (ii)--(iv) for the existence of Pareto efficient solutions to \eqref{problem}. Let $\bar{y}\in f(\Omega)$ be such that the section $Y:= [f(\Omega)]_{\bar y}$ is bounded. Then the closure $\overline Y$ of $Y$ is a nonempty compact set. Fix a vector $\lambda\in\mathrm{int} \, \mathbb{R}_+^m$ and consider the scalar optimization problem
\begin{eqnarray*}
\min_{y\in\overline{Y}}\langle\lambda,y\rangle.
\end{eqnarray*}
which has an optimal solution $\widehat y\in\overline{Y}$.

Assume we have proved that $\Omega\cap f^{-1}(\widehat y)\ne\emptyset$. Take arbitrarily $\widehat x\in\Omega \cap f^{-1}(\widehat y)$ and show that $\widehat x$ is a Pareto efficient solution to problem \eqref{problem}. Arguing by contradiction, suppose that there exists $x\in\Omega$ such that
\begin{eqnarray*}
f(x)&\leqq& f(\widehat x)\quad\textrm{ and }\quad f(x)\ne f(\widehat x).
\end{eqnarray*}
Componentwise it can be equivalently written as
\begin{eqnarray*}
f_i(x)&\le& f_i(\widehat x) \ \textrm{ for } \ i=1,\ldots,n\quad\textrm{ and }\quad f_j(x) \ < \  f_j(\widehat x) \ \textrm{ for some } j.
\end{eqnarray*}
Hence in the case of $f(x)\in Y$ we arrive at the contradiction by
\begin{eqnarray*}
\langle\lambda,f(x)\rangle&<&\langle\lambda,f(\widehat x)\rangle=\langle\lambda,\widehat y\rangle.
\end{eqnarray*}
If otherwise $f(x)\not\in Y,$ we have that $f_i(x)>\bar{y}_i$ for some $i\in \{1,\ldots n\}$, and so
\begin{eqnarray*}
f_i(\widehat x) &=& \widehat y_i \ \le\ \bar{y}_i \ <\ f_i(x)\ \le \ f_i(\widehat x),
\end{eqnarray*}
which is also a contradiction.

It remains to show that the set $\Omega\cap f^{-1}(\widehat y)$ is nonempty provided that either ${K}_{\infty,\leqq \bar{y}}\,(f,\Omega)\subset{K}_{0, \leqq\bar{y}}\,(f,\Omega)$ or ${T}_{\infty,\leqq \bar{y}}\,(f,\Omega)\subset {K}_{0,\leqq\bar{y}}\,(f,\Omega),$ with taking into account that
$\widetilde{K}_{\infty,\leqq\bar{y}}\,(f,\Omega)\subset {K}_{\infty,\leqq\bar{y}}\,(f,\Omega).$ Suppose on the contrary that this claim fails. Denote
\begin{eqnarray*}
X&:=&\Omega \cap f^{-1}(Y) \ = \ \big\{x\in\Omega\,\big|\,f(x)\leqq\bar{y}\big\} \ \ne \ \emptyset,
\end{eqnarray*}
and consider the maximum function $\phi\colon\mathbb{R}^n\rightarrow\mathbb{R}$ defined by \eqref{Eqn3} with its properties mentioned above. Furthermore, it follows from $\widehat y \in\overline{Y}$ that
\begin{eqnarray*}
\inf\limits_{x\in X}\phi(x) &=& 0.
\end{eqnarray*}
There are two cases to be considered.\vspace*{-0.1in}

\subsubsection*{\bf Case~1:} ${K}_{\infty,\leqq\bar{y}}\,(f,\Omega)\subset {K}_{0,\leqq\bar{y}}\,(f,\Omega)$.

By using arguments similar to those employed to establish implication (iii)$\Rightarrow$~(i) of Theorem~\ref{Theorem31} we find a sequence $\{u^k\}\subset X$ satisfying the limiting relationships
\begin{eqnarray*}
\|u^k\|\to\infty,\quad\phi(u^k)\to 0,\quad\textrm{ and }\quad\|u^k\|\nu(u^k)\to 0\quad\textrm{ as }\quad k\to\infty.
\end{eqnarray*}
In particular, $f(u^k)\to\widehat y$ as $k\to\infty$, which yields $\widehat y\in{K}_{\infty,\leqq\bar{y}}\,(f,\Omega).$ This implies, by taking into account the imposed assumption, that $\widehat y\in{K}_{0,\leqq\bar{y}}\,(f,\Omega)$. Thus we arrive at $\widehat y=f(\widehat x)$ for some $\widehat x\in\Omega,$ a contradiction.\vspace*{-0.1in}

\subsubsection*{\bf Case~2:} ${T}_{\infty,\leqq\bar{y}}\,(f,\Omega)\subset {K}_{0,\leqq\bar{y}}\,(f,\Omega)$.

Invoking arguments similar to those used to prove implication (iv)$\Rightarrow$(i) of Theorem~\ref{Theorem31}, we find a sequence $\{v^k\}\subset X \cap\Gamma(f,\Omega)$ satisfying the relationships
\begin{eqnarray*}
\|v^k\|\to\infty\quad\textrm{ and }\quad\phi(v^k)\to 0\quad\textrm{ as }\quad k\to\infty.
\end{eqnarray*}
In particular, we get $f(v^k)\to\widehat y$ as $k\to\infty.$ This gives us by definition that $\widehat y \in {T}_{\infty,\leqq \bar{y}}\,(f,\Omega)$, which yields together with the assumption made that $\widehat y\in{K}_{0,\leqq\bar{y}}\,(f,\Omega)$. Therefore $\widehat y=f(\widehat x)$ for some $\widehat x\in\Omega,$ a contradiction, which completes the proof of the theorem. $\h$\vspace*{0.05in}

In this way we arrive at the verifiable necessary and sufficient conditions for the existence of Pareto efficient solutions in constrained vector optimization with nonsmooth data.

\begin{corollary} \label{Corollary41}
Assume that there is $\bar y\in f(\Omega)$ such that the section $[f(\Omega)]_{\bar y}$ is bounded.
Then the problem~\eqref{problem} admits a Pareto solution provided that one of the following equivalent conditions holds:
\begin{enumerate}
\item[\bf(i)] $f|_{\Omega}$ is proper at the sublevel $\bar y$.
\item[\bf(ii)] $f|_{\Omega}$ satisfies the Palais--Smale condition at the sublevel $\bar y$: $\widetilde{K}_{\infty,\leqq \bar{y}}\,
(f,\Omega)=\emptyset$.
\item [\bf(iii)] $f|_{\Omega}$ satisfies the weak Palais--Smale condition at the sublevel $\bar y$: ${K}_{\infty,\leqq \bar{y}}\,
(f,\Omega)=\emptyset$.
\item [\bf(iv)] $f|_{\Omega}$ is M-tame at the sublevel $\bar{y}$: $T_{\infty,\leqq\bar{y}}\,(f,\Omega)=\emptyset$.
\end{enumerate}
\end{corollary}
{\bf Proof.} This is a direct consequence of Theorems~\ref{Theorem31} and \ref{Theorem41}. $\h$

\subsection{Existence of Geoffrion-properly efficient solutions} The last part of this paper is devoted to the existence of Geoffrion-properly efficient solutions to problem~\eqref{problem}. First we show that the equivalent conditions of Corollary~\ref{Corollary41} do not guarantee the existence of a Geoffrion-properly solution to this problem.

\begin{example}\label{Example41}{\rm Let $\Omega:=\{x\in \mathbb{R}\,|\,-x \leqq 0\}$, and let $f\colon\mathbb{R}\to\mathbb{R}^2$ be a polynomial mapping defined by $f(x)=(f_1(x),f_2(x)):=(-x^2,x)$. It is easy to see that every section of $f(\Omega)$ is bounded and that $f|_{\Omega}$ is proper. Corollary~\ref{Corollary41} tells us that problem \eqref{problem} admits a Pareto efficient solution. In fact, it is not hard to check that the whole set $\Omega$ consists of Pareto efficient solutions to \eqref{problem}.

We claim that the problem \eqref{problem} has no Geoffrion-properly efficient solutions. Indeed, let $\bar x$ be an arbitrary element of $\Omega$. We have to show that for all $M >0$ there exists an index $i\in\{1, 2\}$ and some $x\in\Omega$ with $f_i(x)<f_i(\bar x)$ such that
\begin{eqnarray*}
\frac{f_i(\bar x) - f_i(x)}{f_j(x) - f_j(\bar x)} &>& M
\end{eqnarray*}
whenever $j\in \{1,2\}$ with $f_j(\bar x)<f_j(x)$. To proceed, pick $x>\max\{M, \bar x\}$, $i=1$, and $j=2$. Then we have
$f_1(x) < f_1(\bar x),$ $f_2(\bar x) < f_2(x),$ and
\begin{eqnarray*}
\frac{f_i(\bar x) - f_i(x)}{f_j(x) - f_j(\bar x)} &=& x + \bar x \ > \ M,
\end{eqnarray*}
which verifies the claim.}
\end{example}

As shown in \cite[Theorem~5.2]{Henig82}, a {\em necessary condition} for the existence of Geoffrion-properly efficient solutions to \eqref{problem} is
\begin{eqnarray*}
[f(\Omega)+\mathbb{R}^m_+]^{\oplus}\cap(-\mathbb{R}^m_+)&=&\{0\}.
\end{eqnarray*}
By using \cite[Lemma~3.2.4]{Sawaragi1985}, it can be equivalently rewritten as
\begin{eqnarray}\label{Eqn8}
[f(\Omega)]^{\oplus}\cap(-\mathbb{R}^m_+)&=&\{0\},
\end{eqnarray}
where, given $Y\subset\mathbb{R}^m$, the symbol $[Y]^{\oplus}$ stands for the {\em recession cone} of $Y$ defined by
\begin{eqnarray*}
[Y]^{\oplus}:=\big\{d\in\mathbb{R}^m\,\big|\,\exists\,\{(t_k,y^k)\}\subset\mathbb{R}_+\times Y\;\text{ such that }\;t_k \to 0 \;\textrm{ and }\; t_k y^k\to d\big\}.
\end{eqnarray*}

The next result provides {\em sufficient conditions} for the existence of Geoffrion-properly efficient solutions to the constrained vector optimization problem \eqref{problem}.

\begin{theorem}\label{Theorem42} Under the validity of the necessary optimality condition \eqref{Eqn8}, the following equivalent conditions are sufficient for the existence of a Geoffrion-properly efficient solution to the vector optimization problem \eqref{problem}:
\begin{enumerate}
\item [\bf(i)] $f|_{\Omega}$ is proper;
\item[\bf(ii)] $f|_{\Omega}$ satisfies the Palais--Smale condition;
\item [\bf(iii)] $f|_{\Omega}$ satisfies the weak Palais--Smale condition;
\item [\bf (iv)] $f|_{\Omega}$ is M-tame.
\end{enumerate}
\end{theorem}
{\bf Proof.} Let us first check that conditions (i)--(iv) are indeed equivalent in the setting under consideration. By using Corollary~\ref{Corollary31}, it suffices to show that every section of $f(\Omega)$ is bounded. Supposing the contrary gives us a point $y\in\mathbb{R}^m$ for which the section $[f(\Omega)]_y$ is unbounded. Then we find a sequence $\{x^k\}\subset\Omega$ such that $f(x^k)\leqq y$ whenever $k\in\mathbb{N}$ and $\|f(x^k)\|\to\infty$ as $k\to\infty$. Denoting
\begin{eqnarray*}
t_k=\frac{1}{\|f(x^k)\|}\quad\textrm{ and }\quad d^k := t_k f(x^k),\quad k\in\mathbb{N},
\end{eqnarray*}
we have $t^k\to 0$ as $k\to\infty$  and $\|d^k\|=1$ for all $k\in \mathbb{N}$. Without loss of generality, assume that the sequence $\{d^k\}$ converges to some $d\in\mathbb{R}^m$ with $\|d\|=1$. It follows from the definition that $d\in[f(\Omega)]^{\oplus}$. Since $f(x^k)\leqq y$ for all $k\in\mathbb{N}$, we arrive at $d\in-\mathbb{R}^m_+$, a contradiction.

Now we are ready to prove that the set of Geoffrion-properly efficient solutions to \eqref{problem} is nonempty. Invoking \cite[Theorem~3.2]{Benson1979} and \cite[Theorems~2.1~and~5.1]{Henig82}, it suffices to show that the set $f(\Omega)+\mathbb{R}^m_+$ is closed. To proceed, we deduce from Corollary~\ref{Corollary31} that every section of $f(\Omega)$ is compact. Pick an arbitrary sequence $\{a^k\}\subset f(\Omega)+\mathbb{R}^m_+$ that converges to some $a\in\mathbb{R}^m$ and find sequences $\{y^k\}\subset f(\Omega)$ and $\{d^k\}\subset\mathbb{R}^m_+$ such that $a^k=y^k+d^k$ for all $k\in\mathbb{N}$. Since the sequence $\{a^k\}$ is convergent, there is $\bar a\in\mathbb{R}^m$ with $a^k\leqq\bar a$ for all $k\in\mathbb{N}$. It clearly follows that $y^k\leqq\bar a$ whenever $k\in\mathbb{N}$, and thus $\{y^k\}\subset[f(\Omega)]_{\bar{a}}$. The compactness of $[f(\Omega)]_{\bar{a}}$ gives us a subsequence of $\{y^k\}$, which converges to some $y\in[f(\Omega)]_{\bar{a}}$. This implies that $\{d^k\}$ is also convergent to some $d\in\mathbb{R}^m_+$, and therefore
\begin{eqnarray*}
a&=&y+d\in[f(\Omega)]_{\bar{a}}+\mathbb{R}^m_+\subset f(\Omega)+\mathbb{R}^m_+,
\end{eqnarray*}
which completes the proof of the theorem.$\h$\vspace*{0.05in}

We end this section with the following remarks clarifying relationships of the obtained results with the existence other types of properly efficient solutions to \eqref{problem}.

\begin{remark}{\rm

{\bf(i)} It has been realized in vector optimization (see, e.g., \cite{jahn04}) that in the setting under consideration the concept of Geoffrion-properly efficient solutions agrees with the notions of properly efficient solutions in the senses of Benson \cite{Benson1979} and Henig \cite{Henig82}, and that every Geoffrion-properly efficient solution is also properly efficient in the sense of Borwein \cite{Borwein77}.
Thus Theorem~\ref{Theorem42} also provides sufficient conditions for the existence of properly efficient solutions in the senses of Benson, Henig, and Borwein.

{\bf(ii)} It follows  from \cite[Theorem~5.1]{ha10} that problems \eqref{problem} admits a Henig-properly efficient solution if $\Omega=\mathbb{R}^n$ and the objective mapping $f\colon\mathbb{R}^n\to\mathbb{R}^m$ is bounded from below and satisfies the Palais--Smale condition. Recall that $f$ is {\em bounded from below} if there exists a vector $a\in\mathbb{R}^m$ such that
\begin{eqnarray*}
f(\mathbb{R}^n)&\subset&a+\mathbb{R}^m_+.
\end{eqnarray*}
The above definition readily implies that
\begin{eqnarray*}
\big[f(\mathbb{R}^n)\big]^{\oplus}&\subset&\big[a+\mathbb{R}^m_+\big]^{\oplus}\ =\ \mathbb{R}^m_+,
\end{eqnarray*}
and therefore we get $[f(\mathbb{R}^n)]^{\oplus}\cap(-\mathbb{R}^m_+)=\{0\}$. The converse is true when $m=1$ but fails in general. Thus our Theorem~\ref{Theorem42} essentially improves and extends the result of \cite[Theorem~5.1]{ha10}. To illustrate that we get a proper improvement even for unconstrained problems on $\mathbb{R}$ with polynomial objectives, consider problem \eqref{problem} with $\Omega:=\mathbb{R}$ and the cost mapping $f\colon\mathbb{R}\to\mathbb{R}^2$ defined by $f(x):= (x,x^2)$. It is easy to check that
\begin{eqnarray*}
\big[f(\mathbb{R})\big]^{\oplus}\cap\big(-\mathbb{R}^2_+\big)=\{0\}
\end{eqnarray*}
and  that $f$ satisfies the Palais--Smale condition. Theorem~\ref{Theorem42} tells us that this problem admits a Geoffrion-properly efficient solution. However, the mapping $f$ under consideration is not bounded from below on $\mathbb{R}$, and hence the result of \cite{ha10} is not applicable in this case.}
\end{remark}

\section{Conclusions}

This paper demonstrates that developing a novel approach of variational analysis and generalized differentiation to the existence of global optimal solutions to constrained problems of vector optimization allows us to derive truly new results in this area in both smooth and nonsmooth settings. In this way we show that the developed variational approach leads us to verifiable necessary optimality conditions for the existence of Pareto efficient solutions as well as sufficient conditions for the existence of properly efficient solutions to general constrained problems with locally Lipschitzian cost mappings. In particular, the obtained results dramatically improves the very recent existence theorems of Pareto efficient solutions established in \cite{Kim2018} for unconstrained problems with polynomial cost mappings by using techniques of semialgebraic geometry and polynomial optimization.

We see the following natural directions of future developments of the variational approach to the existence theorems in problems of multiobjective optimization.

1. Avoiding the Lipschitz continuity assumption on cost mappings by considering vector optimization problems with merely continuous and also order semicontinuous cost mappings that frequently arise in applications. According to the scheme implemented above, this requires further investigations of fundamental issues of generalized differential calculus and necessary optimality conditions dealing with non-Lipschitzian mappings and the like.

2. Studying optimization and equilibrium problems with set-valued cost mappings, which are at the core of most recent developments in multiobjective optimization and practical applications to various models in economics, finance, behavioral sciences, etc.; see, e,g., the monographs \cite{khan2015, Mordukhovich2018} and the references therein.

3. Considering vector and set-valued optimization problems in infinite-dimensional spaces. This would open the gate to cover, in particular, various dynamical equilibrium models arising in macroeconomic, mechanics, and systems control governed by constrained evolution equations, inclusions, variational conditions, etc; see, e.g., \cite{Mordukhovich2006,Mordukhovich2018,Xiao2012}. Variation principles and appropriate tools of generalized differentiation provide powerful machinery to successfully proceed in the theoretical developments in this direction with subsequent applications.

\subsection*{Acknowledgments}

The last version of this paper was done while the third author and the fourth author were visiting Department of Applied Mathematics, Pukyong National University, Busan, Republic of Korea in April 2018. These authors would like to thank the department for  hospitality and support during their stay.

\end{document}